\documentclass[10pt]{llncs}
\usepackage{amssymb}
\usepackage{amsmath}
\usepackage{graphicx}
\usepackage{psfrag}
\usepackage{xspace}
\usepackage[noadjust]{cite}
\usepackage[line,color,poly,all]{xy}
\usepackage{verbatim}

\newcommand{\mathscr}{\mathcal }
\newcommand{\vv}{{\mathbf {v}}}

\newcommand{\uu}{{\mathbf {u}}}

\newcommand{\R}{{\mathbb R}}

\newcommand{\C}{{\mathbb C}}
\newcommand{\Z}{{\mathbb Z}}
\newcommand{\Q}{{\mathbb Q}}

\newcommand{\bG}{{\mathbf  G}}

\newcommand{\MMM}{{\mathscr M}}
\newcommand{\HHH}{{\mathfrak H}}

\renewcommand{\O}{{\mathscr O}}

\newcommand{\SL}{{\mathrm{SL}}}
\newcommand{\GL}{{\mathrm{GL}}}

\newcommand{\Res}{{\mathrm{Res}}}

\newcommand{\Gal}{{\mathrm{Gal}}}
\newcommand{\Vor}{Vorono\v{\i}\xspace}

\newcommand{\red}{\text{red}}

\renewcommand{\setminus}{\smallsetminus}
\renewcommand{\epsilon}{\varepsilon}
\newcommand{\cC}{\mathcal{C}}
\newcommand{\ep}{\epsilon}

\newcommand{\bep}{\bar{\ep}}
\newcommand{\s}{\sqrt{2}}

\newcommand{\mat}[1]{\Bigl[\begin{array}{cc} #1 \end{array}\Bigr]}
\newcommand{\matthree}[1]{\biggl[\begin{array}{ccc} #1 \end{array}\biggr]}
\newcommand{\RR}{\mathbb R}
\newcommand{\ZZ}{\mathbb Z}
\newcommand{\QQ}{\mathbb Q}
\newcommand{\CC}{\mathbb C}
\newcommand{\OO}{\mathcal O}
\DeclareMathOperator{\sgn}{sgn}
\DeclareMathOperator{\Mat}{Mat}
\DeclareMathOperator{\Size}{Size}
\DeclareMathOperator{\Norm}{Norm}
\newcommand{\del}{\partial}
\newcommand{\lsp}[1]{{}^{#1}\!}
\newcommand{\bcC}{\bar{\cC}}
\newcommand{\longisomto}{\overset{\sim}{\longrightarrow}}
\begin{document}
\title{Hecke operators and Hilbert modular
forms}
\titlerunning{Hecke operators and Hilbert modular forms}  
%
\author{Paul E. Gunnells and Dan Yasaki}
\authorrunning{Paul E. Gunnells and Dan Yasaki}   
%
\tocauthor{Paul E. Gunnells, Dan Yasaki (UMass Amherst)}
\institute{University of Massachusetts Amherst, Amherst, MA 01003, USA}

\maketitle              

\begin{abstract}
Let $F$ be a real quadratic field with ring of integers $\O$ and with
class number $1$.  Let $\Gamma$ be a congruence subgroup of $\GL_{2}
(\O)$.  We describe a technique to compute the action of the Hecke
operators on the cohomology $H^{3} (\Gamma; \C)$.  For $F$ real
quadratic  this cohomology group contains the cuspidal cohomology
corresponding to cuspidal Hilbert modular forms of parallel weight
$2$.  Hence this technique gives a way to compute the Hecke action on
these Hilbert modular forms.
\end{abstract}

\section{Introduction}
\subsection{}\label{ss:intro}
Let $\bG$ be a reductive algebraic group defined over $\Q $, and let
$\Gamma \subset \bG (\Q )$ be an arithmetic subgroup.  Let $Y = \Gamma
\backslash X$ be the locally symmetric attached to $G = \bG (\R)$ and
$\Gamma$, where $X$ is the global symmetric space, and let $\MMM$ be a
local system on $Y$ attached to a rational finite-dimensional complex
representation of $\Gamma$.  The cohomology $H^{*} (Y; \MMM)$ plays an
important role in number theory, through its connection with
automorphic forms and (mostly conjectural) relationship to
representations of the absolute Galois group $\Gal (\overline{\Q} /\Q
)$ (cf.~\cite{franke,agg,apt,vgt2}).  This relationship is revealed in
part through the action of the \emph{Hecke operators} on the
cohomology spaces.  Hecke operators are endomorphisms induced from a
family of correspondences associated to the pair $(\Gamma ,\bG (\Q
))$; the arithmetic nature of the cohomology is contained in the
eigenvalues of these linear maps.

For $\Gamma \subset \SL_{n} (\Z )$, \emph{modular symbols} provide a
concrete method to compute the Hecke eigenvalues in $H^{\nu } (Y;\MMM
)$, where $\nu = n (n+1)/2-1$ is the top nonvanishing degree
\cite{manin,ash.rudolph}.  Using modular symbols many people have
studied the arithmetic significance of this cohomology group,
especially for $n=2$ and $3$
\cite{cremona,stein,agg,apt,exp.ind,vgt2}; these are the
only two values of $n$ for which $H^{\nu } (Y;\MMM )$ can contain
cuspidal cohomology classes, in other words cohomology classes coming
from cuspidal automorphic forms on $\GL (n)$.  Another setting where
automorphic cohomology has been profitably studied using modular
symbols is that of $\Gamma \subset \SL_{2} (\O)$, where $\O$ is the
ring of integers in a complex quadratic field
\cite{cremona2,crem.whitley,bygott,lingham}.  In this case $Y$ is a
three-dimensional hyperbolic orbifold; modular symbols allow
investigation of $H^{2} (Y;\MMM)$, which again contains cuspidal
cohomology classes.

\subsection{}
Now let $F$ be a real quadratic field with ring of integers $\O$, and
let $\bG$ be the $\Q$-group $\Res_{F/\Q} (\GL_{2})$.  Let $\Gamma
\subseteq \bG (\Q)$ be a congruence subgroup.  In this case we have $X
\simeq \HHH \times \HHH \times \R$, where $\HHH$ is the upper
halfplane (\S\ref{sec:voronoi}).  The locally symmetric space $Y$ is
topologically a circle bundle over a Hilbert modular surface, possibly
with orbifold singularities if $\Gamma$ has torsion.  The cuspidal
cohomology of $Y$ is built from cuspidal Hilbert modular forms.  Hence
an algorithm to compute the Hecke eigenvalues on the cuspidal
cohomology gives a topological technique to compute the Hecke
eigenvalues of such forms.  But in this case there is a big difference
from the setting in \S\ref{ss:intro}: the top degree cohomology occurs
in degree $\nu = 4$, but the cuspidal cohomology appears in degrees
$2,3$.\footnote{The reader is probably more familiar with the case of
$\bG ' = \Res_{F/\Q}\SL_{2}$.  In this case the locally symmetric
space is a Hilbert modular surface, and the cuspidal Hilbert modular
forms contribute to $H^{2}$.  Our symmetric space is slightly larger
since the real rank of $\bG$ is larger than that of $\bG '$.  However,
regardless of whether one studies the Hilbert modular surface or our
$\GL_2$ symmetric space, the cusp forms contribute to the cohomology
in degree one below the top nonvanishing degree.}  Thus modular symbols cannot
``see'' the cuspidal Hilbert modular forms, and cannot directly be
used to compute the Hecke eigenvalues.

\subsection{} In this article we discuss a technique, based on
constructions in \cite{experimental}, that in practice allows one to
compute the Hecke action on the cohomology space $H^{3} (Y; \C)$.
Moreover it is easy to modify our technique to compute with other
local systems; all the geometric complexity occurs for trivial
coefficients.  Here we must stress the phrase \emph{in practice},
since we cannot prove that our technique will actually work.
Nevertheless, the ideas in \cite{experimental} have been successfully
used in practice \cite{computation,AGM-part-2}, and the modifications
presented here have been extensively tested for $F=\Q (\sqrt{2}), \Q
(\sqrt{3})$.

The basic idea is the following.  We first identify a finite
topological model for $H^{3} (Y;\C)$, the \emph{\Vor reduced}
cocycles.  This uses a generalization of \Vor 's reduction theory for
positive definite quadratic forms \cite{koecher,ash.sahc}, which
constructs a $\Gamma$-equivariant tessellation of $X$
(\S\ref{sec:voronoi}).  The Hecke operators do not act directly on
this model, and to accommodate the Hecke translates of reduced
cocycles we work with a larger model for the cohomology, the
(infinite-dimensional) space $S_{1} (\Gamma)$ of \emph{$1$-sharblies}
modulo $\Gamma $ (\S\ref{sec:sharbly}).  The space $S_{1} (\Gamma)$ is
part of a homological complex $S_{*} (\Gamma)$ with Hecke action that
naturally computes the cohomology of $Y$.  Any \Vor reduced cocycle in
$H^{3}$ gives rise to a $1$-sharbly cycle, which allows us to identify
a finite dimensional subspace $S^{\red}_{1} (\Gamma)\subset S_{1}
(\Gamma)$.

The main construction is then to take a general $1$-sharbly cycle
$\xi$ and to modify it by subtracting an appropriate coboundary to
obtain a homologous cycle $\xi '$ that is closer to being \Vor reduced
(\S\ref{sec:reduction}).  By iterating this process, we eventually
obtain a cycle that lies in our finite-dimensional subspace
$S^{\red}_{1} (\Gamma)$.  Unfortunately, we are unable to prove that
at each step the output cycle $\xi '$ is better than the input cycle
$\xi $, in other words that it is somehow ``more reduced.''  However,
in practice this always works.

The passage from $\xi$ to $\xi '$ is based on ideas found in
\cite{experimental}, which describes an algorithm to compute the Hecke
action on $H^{5}$ of congruence subgroups of $\SL_{4} (\Z)$.  The
common feature that this case has with that of subgroups of $\GL_{2}
(\O)$ is that the cuspidal cohomology appears in the degree one less
than the highest.  This means that from our point of view the two
cases are geometrically very similar.  There are some complications,
however, coming from the presence of non-torsion units in $\O$,
complications leading to new phenomena requiring ideas not found in
\cite{experimental}.  This is discussed in \S\ref{sec:comments}.  We
conclude the article by exhibiting the reduction of a $1$-sharbly to a
sum of \Vor reduced $1$-sharblies where the base field is $\Q
(\sqrt{2})$ (\S\ref{sec:example}).

We remark that there is another case sharing these same geometric
features, namely that of subgroups of $\GL_{2} (\O_{K})$, where $K$ is
a complex quartic field.  We are currently applying the algorithm in
joint work with F.~Hajir and D.~Ramakrishnan for $K=\Q (\zeta_{5})$ to
compute the cohomology of congruence subgroups of $\GL_{2} (\OO_{K})$
and to investigate the connections between automorphic cohomology and
elliptic curves over $K$.  Details of these cohomology computations,
including some special features of the field $K$, will appear in
\cite{reduced}; the present paper focuses on the Hilbert modular case.

Finally, we remark that there is a rather different method to compute
the Hecke action on Hilbert modular forms using the Jacquet--Langlands
correspondence.  For details we refer to work of L.~Demb{\'e}l{\'e}
\cite{lassina1,lassina2}.  However, the Jacquet--Langlands technique
works only with the complex cohomology of subgroups of $\GL_{2} (\O)$,
whereas our method in principle allows one to compute with torsion
classes in the cohomology.

\section{Background}
\subsection{}
Let $F$ be a real quadratic field with class
number 1.  Let $\OO \subset F$ denote the ring of integers.  Let $\bG$
be the $\Q$-group $\Res_{F/\QQ}(\GL_2)$ and let $G=\bG(\RR)$ the
corresponding group of real points.  Let $K\subset G$ be a maximal
compact subgroup, and let $A_G$ be the identity component of the
maximal $\QQ$-split torus in the center of $G$.  Then the symmetric
space associated to $G$ is $X=G/KA_G$.  Let $\Gamma \subseteq
\GL_2(\OO)$ be a finite index subgroup. 

In \S\ref{sec:voronoi} we present an explicit model of $X$ in terms
of positive-definite binary quadratic forms over $F$ and construct a
$\GL_2(\OO)$-equivariant tessellation of $X$ following
\cite{koecher,ash.sahc}.  Section~\ref{sec:sharbly} recalls the
sharbly complex \cite{LS,Aunst,experimental}.

\subsection{\Vor{} polyhedron}\label{sec:voronoi}
Let $\iota_1,
\iota_2$ be the two real embeddings of $F$ into $\RR$.  These maps
give an isomorphism $F\otimes_{\QQ} \RR \simeq \RR^2$, and more generally, an
isomorphism
\begin{equation}\label{eq:GL2pre}
G\longisomto \GL_2(\R) \times \GL_2(\R).
\end{equation} 
When the meaning is clear from the context, we use $\iota_1$,
$\iota_2$ to denote all such induced maps.  In particular,
\eqref{eq:GL2pre} is the map
\begin{equation}\label{eq:GL2}
g \longmapsto (\iota_1(g),\iota_2(g)).
\end{equation}
Under this identification, $A_G$ corresponds to $\{(rI,rI)\mid r
>0\}$, where $I$ is the $2\times 2$ identity matrix.

Let $C$ be the cone of real positive definite binary quadratic forms,
viewed as a subset of $V$, the $\RR$-vector space of $2 \times 2$ real
symmetric matrices.  The usual action of $\GL_2(\RR)$ on $C$ is given
by \begin{equation}\label{eq:cone} (g \cdot \phi)(v) = \phi(\lsp{t}gv
), \quad \text{where $g \in \GL_2(\RR)$ and $\phi \in C$.}
\end{equation}
Equivalently, if $A_\phi$ is the symmetric matrix representing $\phi$,
then $g \cdot \phi = g A_\phi \lsp{t}g$.
In particular a coset $gO(2) \in \GL_2(\RR)/O(2)$ can be viewed
as the positive definite quadratic form associated to the symmetric
matrix $g\,\lsp{t}g$.
  
Let $\cC=C\times C$.  Then \eqref{eq:GL2} and \eqref{eq:cone} define
an action of $G$ on $\cC$. Specifically, $g \cdot (\phi_1,\phi_2) =
(\alpha_1,\alpha_2)$, where $\alpha_i$ is represented by $\iota_i(g)
A_{\phi_i}\iota_i(\lsp{t}g)$.  Let $\phi_0$ denote the quadratic form
represented by the identity matrix.  Then the stabilizer in $G$ of
$(\phi_0,\phi_0)$ is a maximal compact subgroup $K$.  The group $A_G$
acts on $\cC$ by positive real homotheties, and we have
\[X=\cC/\R_{>0}=(C\times C)/\RR_{>0} \simeq \HHH \times \HHH \times
\RR,\] where $\HHH$ is the upper halfplane.

Let $\bcC$ denote the closure of $\cC$ in $V \times V$.  Each vector
$w\in \R^2$ gives a rank 1 positive semi-definite form $w \,\lsp{t}w$
(here $w$ is regarded as a column vector).  Combined with $\iota_1$
and $\iota_2$, we get a map $L:\OO^2 \to \bcC$ given by
\begin{equation}
L(v)= \left(\iota_1(v) \cdot \lsp{t}(\iota_1(v)),\iota_2(v) \cdot \lsp{t}(\iota_2(v))\right).
\end{equation}
Let $R(v)$ be the ray $\RR_{>0} \cdot L(v) \subset \bcC$.  Note that 
\[L(cv) = (\iota_1(c)^2L_1(v),\iota_2(c)^2L_2(v))\] 
so that if $c \in \QQ$, then $L(cv)\in R(v)$, and in particular
$L(-v)=L(v)$.  The set of \emph{rational boundary components}
$\cC_{1}$ of $\cC$ is the set of rays of the form $R (v)$, $v\in
F^{2}$ \cite{ash.sahc}.  These are the rays in $\bcC$ that correspond
to the usual cusps of the Hilbert modular variety.

Let $\Lambda \subset V \times V$ be the lattice 
\[
\Lambda=\Bigl\{(\iota_1(A),\iota_2(A))\Bigm| A=\mat{a&c\\c&b},
\quad a,b,c \in \OO \Bigr\}.
\]  
Then $\GL_2(\OO)$ preserves $\Lambda$.
\begin{definition}
The \emph{\Vor{} polyhedron $\Pi$} is the closed convex hull in $\bcC$
of the points $\cC_{1}\cap \Lambda \setminus \{0\}$.
\end{definition}
Since $F$ has class number $1$, one can show that any vertex of $\Pi$
has the form $L (v)$ for $v\in \O^{2}$.  We say that $v\in\O^{2}$ is
\emph{primitive} if $L (v)$ is a vertex of $\Pi$.   Note that $v$
is primitive only if $L (v)$ is primitive in the usual sense as a
lattice point in $\Lambda$.

By construction $\GL_2(\OO)$ acts on $\Pi$.  By taking the cones on
the faces of $\Pi$, one obtains a \emph{$\Gamma$-admissible
decomposition} of $\cC$ for $\Gamma =\GL_{2} (\OO)$ \cite{ash.sahc}.
Essentially this means that the cones form a fan in $\bcC$ and that
there are finitely many cones modulo the action of $\GL_{2} (\OO)$.
Since the action of $\GL_{2} (\OO)$ commutes with the homotheties,
this decomposition descends to a $\GL_{2} (\OO)$-equivariant
tessellation of $X$.\footnote{If one applies this construction to
$F=\Q$, one obtains the Farey tessellation of $\HHH $, with tiles
given by the $\SL_{2} (\Z)$-orbit of the ideal geodesic triangle with
vertices at $0,1,\infty$.  }

We call this decomposition the \emph{\Vor decomposition}.  We call the
cones defined by the faces of $\Pi$ \emph{\Vor{} cones}, and we refer
to the cones corresponding to the facets of $\Pi$ as \emph{top cones}.
The sets $\sigma \cap \cC$, as $\sigma$ ranges over all top cones,
cover $\cC$.  Given a point $\phi \in \cC$, there is a finite
algorithm that computes which \Vor cone contains $\phi$ \cite{Gmod}.

For some explicit examples of the \Vor decomposition over real
quadratic fields, we refer to \cite{Ong} (see also \S\ref{sec:example}).

\subsection{The sharbly complex}\label{sec:sharbly}
Let $S_k$, $k\geq 0$, be the $\Gamma$-module $A_{k}/C_{k}$, where $A_{k}$ is the
set of formal $\ZZ$-linear sums of symbols $[v]=[v_1, \cdots,
v_{k+2}]$, where each $v_i$ is in $F^2$, and $C_{k}$ is the submodule
generated by
\begin{enumerate}
\item $[v_{\sigma(1)}, \cdots, v_{\sigma(k+2)}]-\sgn(\sigma)[v_1,
\cdots, v_{k+2}]$,
\item $[v, v_2, \cdots, v_{k+2}] - [w, v_2, \cdots v_{k+2}]$ if $R (v)
=  R(w)$, and 
\item $[v]$, if $v$ is \emph{degenerate}, i.e., if $v_1, \cdots ,
v_{k+2}$ are contained in a hyperplane.
\end{enumerate}  
We define a boundary map $\partial\colon S_{k+1} \to S_{k}$ by
\begin{equation}\label{eq:boundary}
\del [v_1, \cdots , v_{k+2}] =\sum_{i=1}^{k+2}(-1)^i[v_1, \cdots, \hat{v}_i,\cdots , v_{k+2}].
\end{equation}
This makes $S_{*}$ into a homological complex, called the
\emph{sharbly complex} \cite{ash.sharb}.

The basis elements $\uu =[v_1, \cdots, v_{k+2}]$ are called
\emph{$k$-sharblies}.  Notice that in our class number $1$ setting,
using the relations in $C_k$ one can always find a representative for
$\uu$ with the $v_{i}$ primitive.  In particular, one can always
arrange that each $L(v_i)$ is a vertex of $\Pi$.  When such a
representative is chosen, the $v_i$ are unique up to multiplication by
${\pm 1}$.  In this case the $v_i$---or by abuse of notation the
$L(v_i)$---are called the \emph{spanning vectors} for $\uu$.
\begin{definition}
A sharbly is \emph{\Vor{} reduced} if its spanning vectors are a
subset of the vertices of a \Vor{} cone.
\end{definition}
The geometric meaning of this notion is the following.  Each sharbly
$\uu $ with spanning vectors $v_{i}$ determines a closed cone $\sigma
(\uu)$ in $\bcC$, by taking the cone generated by the points $L
(v_{i})$.  Then $\uu$ is reduced if and only if $\sigma (\uu)$ is
contained in some \Vor cone.  It is clear that there are finitely many
\Vor reduced sharblies modulo $\Gamma$.

Using determinants, we can define a notion of size for $0$-sharblies:
\begin{definition}
Given a $0$-sharbly $\vv$, the \emph{size $\Size (\vv)$ of $\vv$} is
given by the absolute value of the norm determinant of the $2\times 2$
matrix formed by spanning vectors for $\vv$.
\end{definition}
By construction $\Size$ takes values in $\ZZ_{>0}$.  We remark that
the size of a $0$-sharbly $\vv$ is related to whether or not $\vv$ is
\Vor reduced, but that in general there exist \Vor reduced
$0$-sharblies with size $>1$.

The boundary map \eqref{eq:boundary} commutes with the action of
$\Gamma$, and we let $S_*(\Gamma)$ be the homological complex of
coinvariants.  Note that $S_*(\Gamma )$ is infinitely generated as a
$\ZZ \Gamma$-module.  One can show
\begin{equation}\label{eq:iso}
H_{k} ((S_{*}\otimes \CC) (\Gamma))  \longisomto H^{4-k} (\Gamma ; \CC
)
\end{equation}
(cf.~\cite{ash.sharb}), with a similar result holding for cohomology
with nontrivial coefficients.  Moreover, there is a natural action of
the Hecke operators on $S_*(\Gamma)$ (cf.~\cite{experimental}).  Thus
to compute with $H^{3} (\Gamma ; \C)$, which will realize cuspidal
Hilbert modular forms over $F$ of weight $(2,2)$, we work with
$1$-sharbly cycles.  We note that the \Vor reduced sharblies form a
subcomplex of $S_{*} (\Gamma)$ \emph{finitely generated} subcomplex
that also computes the cohomology of $\Gamma$ as in \eqref{eq:iso}.
This is our finite model for the cohomology of $\Gamma$.

\section{The reduction algorithm}\label{sec:reduction}
\subsection{The strategy}\label{sec:strategy}
The general idea behind our algorithm is simple.  To compute the
action of a Hecke operator on the space of $1$-sharbly cycles, it
suffices to describe to an algorithm that writes a general $1$-sharbly
cycle as a sum of \Vor{} reduced $1$-sharblies.  Now any basis
$1$-sharbly $\uu$ contains three sub-$0$-sharblies (the \emph{edges}
of $\uu$), and the \Vor{} reduced $1$-sharblies tend to have edges of
small size.  Thus our first goal is to systematically replace all the
$1$-sharblies in a cycle with edges of large size with $1$-sharblies
having smaller size edges.  This uses a variation of the classical
modular symbol algorithm, although no continued fractions are
involved.  Eventually we produce a sum of $1$-sharblies with all edges
\Vor reduced.  However, having all three edges \Vor reduced is
(unfortunately) not a sufficient condition for a $1$-sharbly to be
\Vor reduced.\footnote{This is quite different from what happens with
classical modular symbols, and reflects the infinite units in $\O$.}
Thus a different approach must be taken for such $1$-sharblies to
finally make the cycle \Vor reduced. This is discussed further in
\S\ref{sec:comments}.

\subsection{Lifts}
We begin by describing one technique to encode a $1$-sharbly cycle
using some mild extra data, namely that of a choice of \emph{lifts}
for its edges:

\begin{definition}[\cite{experimental}]
A $2\times 2 $ matrix $M$ with coefficients of $F$ with columns $A_1$,
$A_2$ is said to be a \emph{lift} of a $0$-sharbly $[u,v]$ if
$\{R(A_1),R(A_2)\}=\{R(u),R(v)\}$.
\end{definition}

\begin{figure}
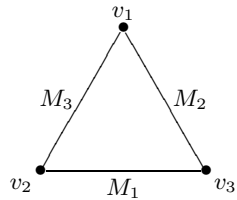

\[{
\xy /r0.5in/:="a",
{\xypolygon3"A"{~>{-}{\bullet}}},
"A1" *+!D{v_1},"A2" *+!UR{v_2},"A3" *+!UL{v_3},
"A1";"A2" **@{} ?*{}="c1" *+!R{M_3},
"A2";"A3" **@{} ?*{}="c2" *+!U{M_1},
"A3";"A1" **@{} ?*{}="c3" *+!L{M_2},
\endxy}
\]
\caption{A $1$-sharbly with lifts}
\label{fig:sharbly}
\end{figure}

The idea behind the use of lifts is the following.  Suppose a linear
combination of 1-sharblies $\xi = \sum a (\uu)\uu\in S_{1}$ becomes a
cycle in $S_{1} (\Gamma)$.  Then its boundary must vanish modulo
$\Gamma$.  In the following algorithm, we attempt to pass from $\xi$
to a ``more reduced'' sharbly $\xi '$ by modifying the edges of each
$\uu$ in the support of $\xi$.  To guarantee that $\xi '$ is a cycle
modulo $\Gamma$, we must make various choices in the course of the
reduction $\Gamma$-equivariantly across the boundary of $\xi$.  This
can be done by first choosing $2\times 2$ integral matrices for each
sub-$0$-sharbly of $\xi$.  We refer to \cite{experimental} for more
details and discussion.  For the present exposition, we merely remark
that we always view a $1$-sharbly $\uu = [v_{1}, v_{2}, v_{3}]$ as a
triangle with vertices labelled by the $v_{i}$ and with a given
(fixed) choice of lifts for each edge (Figure~\ref{fig:sharbly}).  If
two edges $\vv , \vv '$ satisfy $\gamma \cdot \vv = \vv '$, then we
choose the corresponding lifts to satisfy $\gamma M = M'$.  The point
is that we can then work individually with $1$-sharblies enriched with
lifts; we don't have to know explicitly the matrices in $\Gamma$ that
glue the $1$-sharblies into a cycle modulo $\Gamma$.

We emphasize that the lift matrices for any given $1$-sharbly in the
support of $\xi$ are essentially forced on us by the requirement that
$\xi$ be a cycle modulo $\Gamma$.   There is almost no
flexibility in choosing them.  Such matrices form an essential part of
the input data for our algorithm.

\subsection{Reducing points}

\begin{definition}\label{def:reducingpoint}
Let $\vv$ be a $0$-sharbly with spanning vectors $\{x,y\}$.  Assume
$\vv$ is not \Vor{} reduced. Then $u\in \OO^{2}\smallsetminus \{0 \}$
is a \emph{reducing point for $\vv$} if the following hold:
\begin{enumerate}
\item $R ( u) \neq R (x), R (y)$.
\item $L(u)$ is a vertex of the unique \Vor{} cone $\sigma$ (not necessarily
top-dimensional) containing the ray $R (x+y)$.
\item If $x=ty$ for some $t\in F^\times$, then $u$ is in the span of $x$. 
\item Of the vertices of $\sigma$, the point $u$ minimizes the sum of
the sizes of the $0$-sharblies $[x,u]$ and $[u,y]$.
\end{enumerate}
\end{definition}
Given a non-\Vor reduced $0$-sharbly $\vv = [x,y]$ and a reducing
point $u$, we apply the relation
\begin{equation}\label{eq:relation}
[x,y] = [x,u]+[u,y]
\end{equation}
in the hopes that the two new $0$-sharblies created are closer to
being \Vor reduced.  Note that choosing $u$ uses the geometry of the
\Vor decomposition instead of (a variation of) the continued fraction
algorithms of \cite{manin,cremona,ash.rudolph}.  Unfortunately we
cannot guarantee that the new $0$-sharblies on the right of
\eqref{eq:relation} are better than $\vv$, but this is true in
practice.

\subsection{$\Gamma$-invariance} \label{sec:normal} The reduction
algorithm proceeds by picking reducing points for non-\Vor{} reduced
edges.  We want to make sure that this is done $\Gamma$-equivariantly;
in other words that if two edges $\vv$, $\vv '$ satisfy $\gamma \cdot
\vv = \vv'$, then if we choose $u$ for $\vv$ we want to make sure that we
choose $\gamma u$ for $\vv '$.

We achieve this by making sure that the choice of reducing point for
$\vv$ only depends on the lift matrix $M$ that labels $\vv $.  The
matrix is first put into \emph{normal form}, which is a unique
representative $M_0$ of the coset $\GL_2(\OO) \backslash M$.  This is
an analogue of Hermite normal form that incorporates the action of the
units of $\OO$.  There is a unique $0$-sharbly associated to $M_0$; We
choose a reducing point $u$ for this $0$-sharbly and translate it back
to obtain a reducing point for $\vv$.  Note that $u$ need not be
unique.  However we can always make sure that the same $u$ is chosen
any time a given normal form $M_{0}$ is encountered, for instance by
choosing representatives of the \Vor cones modulo $\GL_{2} (\OO)$ and
then fixing an ordering of their vertices.

We now describe how $M_{0}$ is constructed from $M$.  Let $\Omega_*$
be a fundamental domain for the action of $(\OO^\times,\cdot)$ on $F^{\times }$.
For $t \in \OO$, let $\Omega_+(t)$ be a fundamental domain for the
action of $(t\OO,+)$ on $F$.
\begin{definition}
A nonzero matrix $M \in \Mat_2(F)$ is in \emph{normal form} if $M$ has one of
the following forms:
\begin{enumerate}
\item $\mat{0&b\\0&0}$, where $b \in \Omega_*$.
\item $\mat{a&b\\0&0}$, where $a \in \Omega_*$ and $b\in F$.
\item $\mat{a&b\\0&d}$, where $a,d \in \Omega_*$ and $b\in \Omega_+(d)$.
\end{enumerate}
\end{definition}
It is easy to check that the normal form for $M$ is uniquely
determined in the coset $\GL_{2} (\OO)\cdot M$.
To explicitly put $M=\mat{a&b\\c&d}$ in normal form, the first step is
to find $\gamma \in \GL_2(\OO)$ such that $\gamma \cdot M$ is upper
triangular.  Such a $\gamma$ can be found after finite computations as
follows.  Let $N\colon F \to \R$ be defined by
$N(\alpha)=|\Norm_{F/\QQ}(\alpha)|$.  If \[0<N(c) < N(a),\] then let
$\alpha \in \OO$ be an element of smallest distance from $a/c$.  Let
\[\gamma'=\mat{0&1\\1&0}\mat{1&-\alpha\\0&1}.\] Then $\gamma' \in
\GL_2(\OO)$ and $\gamma' M=\mat{a'&b'\\c'&d'}$ with $N(c') < N(c)$ and
$N(a')<N(a)$.  Repeating this procedure will yield the desired result.

After a reducing point is selected for $\vv$ and the relation
\eqref{eq:relation} is applied, we must choose lifts for the
$0$-sharblies on the right of \eqref{eq:relation}.  This we do as
follows:

\begin{definition}
Let $[v_1,v_2]$ be a non-reduced $0$-sharbly with lift matrix $M$ and
reducing point $u$.  Then the \emph{inherited lift} $\hat{M}_i$ for
$[v_i,u]$ is the matrix obtained from $M$ by keeping the column
corresponding to $v_i$ and replacing the other column by $u$.
\end{definition}

\subsection{The algorithm} Let $T=[v_1,v_2,v_3]$ be a non-degenerate
sharbly.  Let $M_i$ be the lifts of the edges of $T$ as shown in
Figure~\ref{fig:sharbly}.  The method of subdividing the interior
depends on the number of edges that are \Vor{} reduced.  After each
subdivision, lift data is attached using inherited lifts for the
exterior edges. The lift for each interior edge can be chosen
arbitrarily as long as the same choice is made for the edge to which
it is glued.  We note that steps (I), (II), and (III.1) already appear
in \cite{experimental}, but (III.2) and (IV) are new subdivisions
needed to deal with the complications of the units of $\OO$.
 
\subsubsection{(I) Three non-reduced edges.} If none of the edges are
\Vor{} reduced, then we split each edge by choosing reducing points
$u_1,u_2$, and $u_3$.  In addition, form three additional edges
$[u_1,u_2],[u_2,u_3]$, and $[u_3,u_1]$.  We then replace $T$ by the
four $1$-sharblies
\[[v_1,v_2,v_3]\longmapsto[v_1,u_3,u_2]+[u_3,v_2,u_1]+[u_2,u_1,v_3]+[u_1,u_2,u_3].\]
\subsubsection{(II) Two non-reduced edges.}  If only one edge is \Vor{}
reduced, then we split the other two edges by choosing reducing points
$u_1$ and $u_3$.  We form two additional edges $[u_1,u_3]$ and $\ell$,
where $\ell$ is taken to be either $[v_1,u_1]$ or $[v_3,u_3]$,
whichever has smaller size.  More precisely:
\begin{enumerate}
\item If $\Size([v_1,u_1])\le \Size([u_3,v_3])$, then we form two additional edges $[u_1,u_3]$ and $[v_1,u_1]$, and replace $T$ by the three $1$-sharblies
\[[v_1,v_2,v_3] \longmapsto [v_1,u_3,u_1]+[u_3,v_2,u_1]+[v_1,u_1,v_3].\]
\item Otherwise, we form two additional edges $[u_1,u_3]$ and $[v_3,u_3]$, and replace $T$ by the three $1$-sharblies
\[[v_1,v_2,v_3] \longmapsto [v_1,u_3,v_3]+[u_3,v_2,u_1]+[u_3,u_1,v_3].\]
\end{enumerate}
\subsubsection{(III) One non-reduced edge.} If two edges are
\Vor{} reduced, then we split the other edge by choosing a reducing
point $u_1$.  The next step depends on the configuration of
$\{v_1,v_2,v_3,u_1\}$.
\begin{enumerate}
\item If $[v_2,u_1]$ or $[u_1,v_3]$ is not \Vor{} reduced or
$v_2=tv_1$ for some $v \in F$, then we form one additional edge
$[v_1,u_1]$ and replace $T$ by the two $1$-sharblies
\[[v_1,v_2,v_3]\longmapsto[v_1,v_2,u_1]+[v_1,u_1,v_3].\]
\item Otherwise, a central point $w$ is chosen.  The central point $w$
is chosen from the vertices of the top cone containing the barycenter
of $[v_1,v_2,v_3,w]$ so that it maximizes the number of \Vor{} reduced
edges in the set \[S=\{[v_1,w],[v_2,w],[v_3,w],[u_1,w]\}.\] We do not
allow $v_1,v_2$ or $v_3$ to be chosen as a central point.  We form
four additional edges $[v_1,w],[v_2,w],[u_1,w],$ and $[v_3,w]$ and
replace $T$ by the four $1$-sharblies
\[[v_1,v_2,v_3]\longmapsto[v_1,v_2,w]+[w,v_2,u_1]+[w,u_1,v_3]+[w,v_3,v_1].\]
\end{enumerate}
\subsubsection{(IV) All edges \Vor reduced.}  If all three edges are
\Vor{} reduced, but $T$ is not \Vor{} reduced, then a central point
$w$ is chosen.  The central point $w$ is chosen from the vertices of
the top cone containing the barycenter of $[v_1,v_2,v_3]$ so that it
maximizes the sum $\#E+\#P$, where $E$ is the set of \Vor{} reduced
edges in $\{[v_1,w],[v_2,w],[v_3,w]\}$ and $P$ is the set of \Vor{}
reduced triangles in $\{[v_1,v_2,w],[v_2,v_3,w],[v_3,v_1,w]\}$.  We do
not allow $v_1,v_2$ or $v_3$ to be chosen as a central point.  We form
three additional edges $[v_1,w],[v_2,w],$ and $[v_3,w]$ and replace
$T$ by the three $1$-sharblies
\[[v_1,v_2,v_3]\longmapsto[v_1,v_2,w]+[w,v_2,v_3]+[w,v_3,v_1].\]
\section{Comments}\label{sec:comments}
First, we emphasize that the reducing point $u$ of
Definition~\ref{def:reducingpoint} works in practice to shrink the
size of a 0-sharbly $\vv $, but we have no proof that it will do so.
The difficulty is that Definition~\ref{def:reducingpoint} chooses $u$
using the geometry of the \Vor polyhedron $\Pi $ and not the size of
$\vv$ directly.  Moreover, our experience with examples shows that
this use of the structure of $\Pi$ is essential to reduce the original
$1$-sharbly cycle (cf.~\S\ref{sec:example2}).

Next, as mentioned in \S\ref{sec:strategy}, case (IV) is necessary:
there are $1$-sharblies $T$ with all three edges \Vor reduced, yet $T$
is itself not \Vor reduced.  An example is given in the next section.
The point is that in $\bcC $ the points $L (v)$ and $L (\epsilon v)$
are different if $\epsilon $ is not a torsion unit, but after passing
to the Hilbert modular surface $L (v)$ and $L (\epsilon v)$ define the
same cusp.  This means one can take a geodesic triangle $\Delta$ in
the Hilbert modular surface with vertices at three cusps that by any
measure should be considered reduced, and can lift $\Delta $ to a
$3$-cone in the $\GL_{2}$-symmetric space that is far from being \Vor
reduced.

Finally, the reduction algorithm can be viewed as a two stage process.
When a 1-sharbly $T$ has 2 or 3 non-reduced edges or 1 non-reduced
edge and satisfies the criteria for case 1, then in some sense $T$ is
``far'' from being \Vor{} reduced.  One tries to replace $T$ by a sum
of 1-sharblies that are more reduced in that the edges have smaller
size.  However, this process will not terminate in \Vor{} reduced
sharblies.  In particular, if $T$ is ``close'' to being \Vor{}
reduced, then one must use the geometry of the \Vor{} cones more
heavily.  This is why we need the extra central point $w$ in (III.2)
and (IV).

For instance, suppose $T=[v_1,v_2,v_3]$ is a $1$-sharbly with 1
non-reduced edge such that the criteria for (III.2) are satisfied when
the reducing point is chosen.  One can view choosing the central point
and doing the additional split as first moving the bad edge to the
interior of the triangle, where the splitting no longer needs to be
$\Gamma$-invariant.  The additional freedom allows one to make a
better choice.  Indeed, without the central point chosen wisely, this
does lead to some problems. In particular, there are examples where
$[v_1,u_1]$ is not \Vor{} reduced, and the choice of the reducing
point for this edge is $v_2$, leading to a repeating behavior.  Thus
the distinction had to be made.

\section{The case $F=\QQ (\s)$}\label{sec:example}
\subsection{}
Let $F=\QQ(\s)$ and let $\epsilon=1+\s$, a fundamental unit of norm $-1$.
Computations of H.~Ong \cite[Theorem~4.1.1]{Ong} with positive definite binary quadratic forms
over $F$ allow us to describe the \Vor polyhedron
$\Pi$ and thus the \Vor decomposition of $\cC$:

\begin{proposition}[{\cite[Theorem~4.1.1]{Ong}}]
Modulo the action of $\GL_{2} (\OO)$, there are two inequivalent top
\Vor cones.  The corresponding facets of $\Pi$ have $6$ and $12$
vertices, respectively.
\end{proposition}

We fix once and for all representative $6$-dimensional cones $A_0$ and
$A_1$.  To describe these cones, we give sets of points $S\subset
\OO^{2}$ such that the points $\{L (v)\mid v\in S \}$ are the vertices
of the corresponding face of $\Pi $. Let $e_{1}, e_{2}$ be the
canonical basis of $\OO^{2}$.  Then we can take $A_{0}$ to
correspond to the $6$ points
\[
e_{1}, e_{2}, e_{1}-e_{2},\bep e_{1}, \bep e_{2}, \bep (e_{1}-e_{2}),
\]
and $A_{0}$ to correspond to the $12$ points
\[
e_{1}, e_{2}, \bep e_{1}, \bep e_{2}, e_{1}-e_{2}, e_{1}+\bep e_{2},
e_{2}+\bep e_{1}, \bep (e_{1}+e_{2}), \alpha, \beta, \bep \alpha, \bep \beta,
\]
where $\alpha =e_{1}-\s e_{2}$, $\beta =e_{2}-\s e_{1}$.
Since $A_{1}$ is not a simplicial cone, there exist basis sharblies
that are \Vor reduced but do not correspond to \Vor cones.

Now we consider cones of lower dimension.  Modulo $\GL_2(\OO)$, every 
2-dimensional \Vor{} cone either lies in $\bcC \setminus \cC$ or is
equivalent to the cone corresponding to $\{e_{1}, e_{2}\}$.  The $\GL_2(\OO)$-orbits of 3-dimensional \Vor{}
cones are represented by $\{e_{1}, e_{2} \}\cup U$, where $U$ ranges over
\[
\{e_{1}-e_{2} \}, \{\bep e_{1}\},
\{\bep (e_{1}-e_{2}) \}, \{e_{1}-\s e_{2}, e_{2}-\s e_{1}\}, \{e_{1}+\bep e_{2} \}.
\]

Note that all but one of the $3$-cones are simplicial.

\subsection{}\label{sec:example2}
Now we consider reducing a $1$-sharbly $T$.
Let us represent $T$ by a $2\times 3$ matrix whose columns
are the spanning vectors of $T$.  We take $T$ to be 
\[T=\matthree{   \s + 3 & 4\s + 4&  3\s - 4\\
       \s & 5\s - 1& -3\s - 5},\] and we choose arbitrary initial
lifts for the edges of $T$.  This data is typical of what one
encounters when trying to reduce a $1$-sharbly cycle modulo $\Gamma $.

The input $1$-sharbly $T$ has 3 non-reduced edges with edge sizes given by the vector
$[5299,529,199]$.  The first pass of the algorithm follows (I) and
splits all 3 edges, replacing $T$ by the sum $S_1+S_2+S_3+S_4$,
where
\[S_1=
    \matthree{ \s + 3 &-\s - 1    &  1\\
          \s &    -\s    &  0},
  \quad S_2= \matthree{4\s + 4 &      0  &-\s - 1\\
    5\s - 1 & -\s - 1 &     -\s},\]
\[    S_3=\matthree{ 3\s - 4 &       1 &       0\\
    -3\s - 5 &       0 &  -\s - 1},
   \quad S_4=\matthree{     0   &   1 &-\s - 1\\
    -\s - 1 & 0& -\s} .
\] 
We compute that $\Size(S_1)=[2,2,8]$, $\Size(S_2)=[1,1,16]$,
$\Size(S_3)=[1,2,7]$, and $\Size(S_4)=[2,1,1]$. Notice that the
algorithm replaces $T$ by a sum of sharblies with edges of
significantly smaller size.  This kind of performance is typical, and
looks similar to the performance of the usual continued fraction
algorithm over $\Z$.  Note also that $S_{4}$, which is the $1$-sharbly
spanned by the three reducing points of the edges $T$, also has edges
of very small size.  This reflects our use of
Definition~\ref{def:reducingpoint} to choose the reducing points;
choosing them without using the geometry of $\Pi$ often leads to bad
performance in the construction of this $1$-sharbly.

Now $S_4$ has 3 \Vor{} reduced edges, but is itself not \Vor{}
reduced.  The algorithm follows (IV), replaces $S_4$ by $R_1+R_2+R_3$,
and now each $R_i$ is \Vor{} reduced.  

The remaining $1$-sharblies $S_1$, $S_2$, and $S_3$ have only 1
non-reduced edge.  They are almost reduced in the sense that they
satisfy the criteria for (III.2).  The algorithm replaces
$S_1$ by $O_1+O_2+O_3+O_4$, where $O_1$ and $O_2$ are degenerate and
$O_3$ and $O_4$ are \Vor{} reduced.  The $1$-sharbly $S_2$ is replaced
by a $P_1+P_2+P_3+P_4$, and each $P_i$ is \Vor{} reduced.  $S_3$ is
replaced by $Q_1+Q_2+Q_3+Q_4$, where $Q_1$ and $Q_2$ are degenerate,
$Q_3$ is \Vor{} reduced, and $Q_4$ is not \Vor{} reduced.
This $1$-sharbly is given by 
\[Q_4=\matthree{  -\s + 1 &       0&  3\s - 4\\
         2\s + 3 & -\s - 1 &-3\s - 5}\] and has 3 \Vor{} reduced
edges.  Once again the algorithm is in case (IV), and replaces $Q_4$
by a sum $N_1+N_2+N_3$ of \Vor{} reduced sharblies.

To summarize, the final output of the reduction algorithm applied to $T$ is a sum 
\[N_1+N_2+N_3+O_3+O_4+P_1+P_2+P_3+P_4+Q_3+R_1+R_2+R_3, \quad \text{where}\]
\begin{gather*}
N_1=\matthree{  -\s + 1 &       0   &     0\\
     2\s + 3  & -\s - 1 &-2\s - 3},\\
 N_2=    \matthree{       0 & 3\s - 4   &     0\\
     -\s - 1& -3\s - 5 &-2\s - 3}, \quad\\
N_3  =  \matthree{ 3\s - 4  & -\s + 1  &      0\\
   -3\s - 5  &2\s + 3 &-2\s - 3}, \quad
O_3=\matthree{-\s - 1& -\s - 1  &    1\\
            -1 &    -\s   &   0}, \\
O_4 =       \matthree{-\s - 1 &     1 & \s + 3\\
            -1  &    0  &    \s}, \quad
P_1=\matthree{ 2\s + 3 & 4\s + 4 &       1\\
           \s + 2 & 5\s - 1 &-2\s + 2},\\
 P_2=        \matthree{ 2\s + 3 &       1   &     0\\
           \s + 2& -2\s + 2  & -\s - 1},\quad
P_3=        \matthree{ 2\s + 3  &     0  &-\s - 1\\
         \s + 2 & -\s - 1 &     -\s}, \\
P_4=        \matthree{ 2\s + 3 & -\s - 1 &4\s + 4\\
          \s + 2 &     -\s &5\s - 1}, \quad
Q_3= \matthree{ -\s + 1&       1     &  0\\
        2\s + 3  &     0  &-\s - 1},\\
 R_1=\matthree{  -\s + 1 &       0   &     0\\
     2\s + 3  & -\s - 1 &-2\s - 3},\\
R_2=\matthree{       0 & 3\s - 4   &     0\\
      -\s - 1& -3\s - 5 &-2\s - 3}, \quad \text{and}\\
 R_3   =\matthree{ 3\s - 4  & -\s + 1  &      0\\
    -3\s - 5  &2\s + 3 &-2\s - 3},
\end{gather*}
 and each of the above is \Vor{} reduced.  Some of these $1$-sharblies
correspond to \Vor cones and some don't.  In particular, one can check
that the spanning vectors for $P_3$, $P_4$, $R_1$, and $N_1$ do form
\Vor{} cones, and all others don't.  However, the spanning vectors of
$O_3$ and $O_4$ almost do, in the sense that they are subsets of
3-dimensional \Vor{} cones with four vertices.  

\providecommand{\bysame}{\leavevmode\hbox to3em{\hrulefill}\thinspace}
\providecommand{\MR}{\relax\ifhmode\unskip\space\fi MR }
\providecommand{\MRhref}[2]{%
  \href{http://www.ams.org/mathscinet-getitem?mr=#1}{#2}
}
\providecommand{\href}[2]{#2}

\end{document}